\RequirePackage{fix-cm}
\documentclass{svjour3}                     
\smartqed  
\usepackage{amsmath,amsfonts,amssymb,amsbsy}
\usepackage{natbib}
\usepackage[cp1251]{inputenc}
\usepackage[english]{babel}
\usepackage{graphicx}
\begin{document}

\title{Sinusoidally-driven unconfined compression test for a biphasic tissue}
\author{I.~Argatov}

\institute{I.~Argatov \at
              Institute of Mathematics and Physics, Aberystwyth University, Ceredigion SY23 3BZ, Wales, UK \\
              Tel.: +44-1970-622768\\
              Fax: +44-1970-622826\\
              \email{iva1@aber.ac.uk}
}

\date{Received: date / Accepted: date}

\maketitle

\begin{abstract}
In recent years, a number of experimental studies have been conducted to investigate the mechanical behavior of water-saturated biological tissues like articular cartilage under dynamic loading. 
For in vivo measurements of tissue viability, the indentation tests with the half-sinusoidal loading history were proposed. 
In the present paper, the sinusoidally-driven compression test utilizing either the load-controlled or displace\-ment-controlled loading protocol are considered in the framework of linear biphasic layer model.
Closed-form analytical solutions for the integral characteristics of the compression test are obtained. 

\end{abstract}

\setcounter{equation}{0}

\section{Introduction} 
\label{1biSectionI}

It is well known that the long-term creep and relaxation tests, typically used for determining viscoelastic and  biphasic/poroelastic properties, are not appropriate for rapidly assessing the dynamic biomechanical properties of biological tissues like articular cartilage. For in vivo measurements of tissue viability, \citet{Appleyard_et_al_2001} developed a dynamic indentation instrument, which  employs a single-frequency (20~Hz) sinusoidal oscillatory waveform superimposed on a carrier load. It is to note that the oscillation test requires some time period to elapse before recording the measurement to minimize the influence of the initial conditions. That is why the indentation or compression tests with the half-sinusoidal loading history are so tempting. In particular, the half-sinusoidal history is useful for developing indentation-type scanning probe-based nanodevices reported by \citet{Stolz_et_al_2007}. Also, as a first approximation, such a indentation history can be used for modeling impact tests \citep{ButcherSegalman2000}.

In the present paper, it is assumed that the mechanical response of a time-dependent material can be described in the framework of biphasic layer model. Following \citet{Argatov2012}, we consider sinusoidally-driven flat-ended compression test utilizing either the load-controlled or displacement-controlled loading protocol. Closed-form analytical solutions for the integral characteristics of the compression test are obtained. 

\section{Unconfined compression of a cylindrical biphasic sample}
\label{1biSection2}

\begin{figure}[h!]
    \centering
    \hbox{
    \includegraphics[scale=0.32]{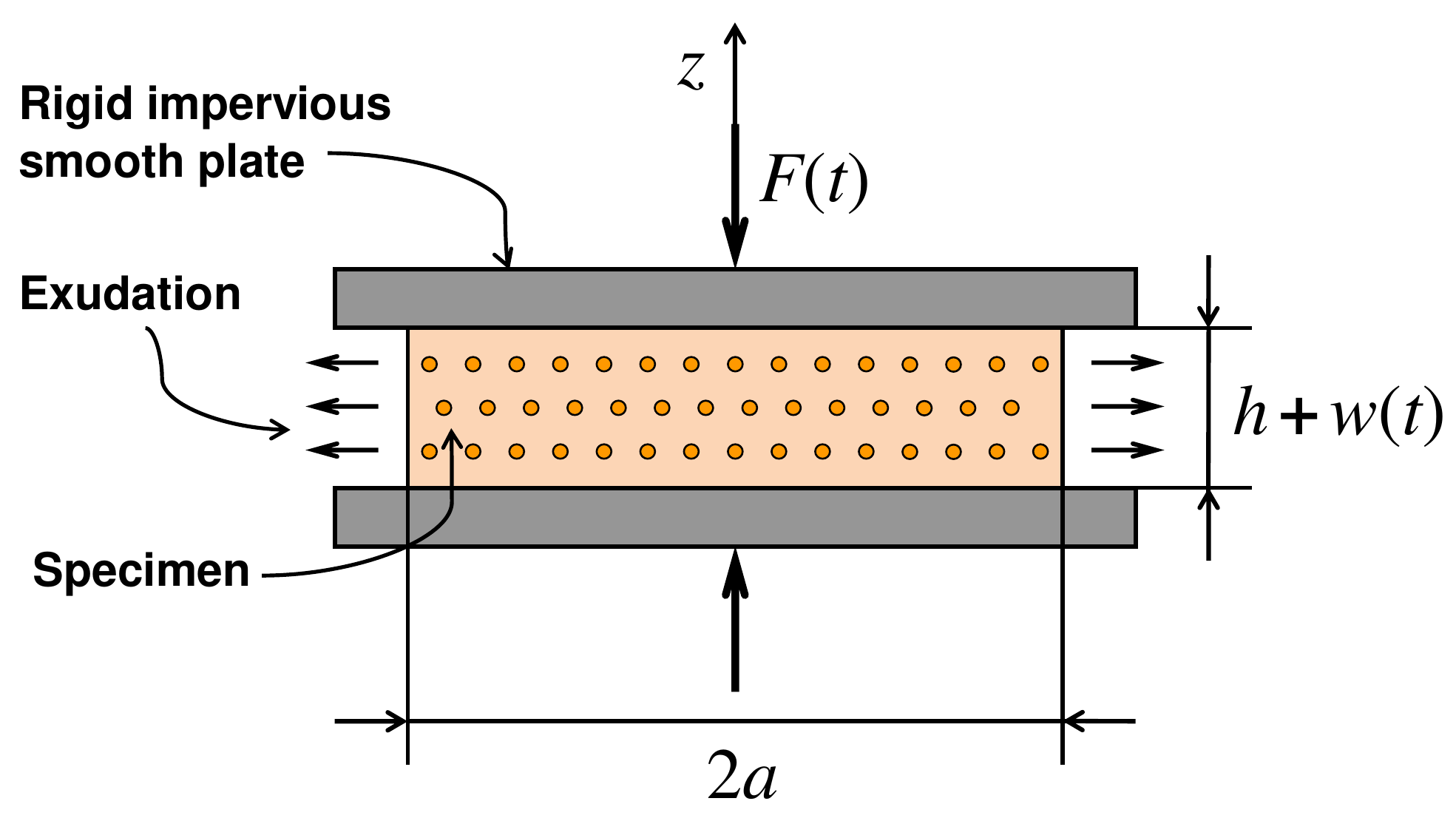}
    }
    \caption{Schematics of the unconfined compression test.    }
    \label{Figure_bi1o0}
\end{figure}

We assume that the unconfined compression test for a soft biological tissue sample can be described by a one-dimensional axisymmetric mathematical model of \citet{Armstrong_et_al_1984} developed in the framework of the linear biphasic theory \citep{Mow_et_al_1980}. In particular, it is assumed that a thin cylindrical specimen is squeezed between two perfectly smooth, impermeable rigid plates such that compressive strain in the axial direction (see, Fig.~\ref{Figure_bi1o0}) is homogeneous. 

Let $F(t)$ and $w(t)$ denote, respectively, the variable compressive force and the vertical displacement of the upper plate with respect to the lower plate ($t$ is the time variable). Then, the compressive strain is defined as 
\begin{equation}
\epsilon(t)=\frac{w(t)}{h},
\label{1bi(2.1)}
\end{equation}
where $h$ is the sample thickness. (Note that compressive strain is negative.)

Let us also introduce the non-dimensional variables 
\begin{equation}
\hat{F}=\frac{F}{\pi a^2\mu_s}, \quad
\hat{t}=\frac{H_A kt}{a^2},
\label{1bi(2.2)}
\end{equation}
where $H_A=\lambda_s+2\mu_s$, $\mu_s$ and $\lambda_s$ are the confined compression modulus and Lam\'e constants of the elastic solid matrix, respectively, $k$ is the tissue permeability, $a$ is the radius of specimen, $\hat{t}$ is the dimensionless time, $\hat{F}$ is the dimensionless force. 

Finally, let $\bar{F}(s)$ and $\bar{\epsilon}(s)$ denote the Laplace transforms with respect to the dimensionless time, and $s$ is the Laplace transform parameter. According to \citet{Armstrong_et_al_1984}, the following relationship holds true:
\begin{equation}
\bar{F}(s)=s\bar{\epsilon}(s)\bar{K}(s),
\label{1bi(2.3)}
\end{equation}
\begin{equation}
\bar{K}(s)=\frac{\displaystyle 3 I_0(\sqrt{s})-\frac{8\mu_s}{H_A}\frac{I_1(\sqrt{s})}{\sqrt{s}}
}{\displaystyle s\Bigl(I_0(\sqrt{s})-\frac{2\mu_s}{H_A}\frac{I_1(\sqrt{s})}{\sqrt{s}}\Bigr)}.
\label{1bi(2.4)}
\end{equation}
Here, $I_0$ and $I_1$ are modified Bessel functions of the first kind. 

By applying the convolution theorem to Eq.~(\ref{1bi(2.3)}), we obtain
\begin{equation}
\hat{F}(\hat{t})=\int\limits_{0-}^{\hat{t}} \frac{d\epsilon(\tau)}{d\tau}\hat{K}(\hat{t}-\tau)\,d\tau,
\label{1bi(2.5)}
\end{equation}
where $\hat{K}(\hat{t})=\mathcal{L}^{-1}\{\bar{K}(s)\}$ is the original function of $\bar{K}(s)$, $t=0-$ is the time moment just preceding the initial moment of loading. In deriving Eq.~(\ref{1bi(2.5)}), we used the formula
$$
\mathcal{L}^{-1}\{s\bar{\epsilon}(s)\}=\frac{d\epsilon(\hat{t})}{d\hat{t}}+\epsilon(0+)\delta(\hat{t}),
$$
where $\delta(\hat{t})$ is the Dirac function. 

By using the results obtained by \citet{Armstrong_et_al_1984}, we will have 
\begin{equation}
\hat{K}(\hat{t})=2(1+\nu_s)+\sum_{n=1}^\infty 2(1+\nu_s)A_n\exp(-\alpha_n^2\hat{t}),
\label{1bi(2.9)}
\end{equation}
\begin{equation}
A_n=\frac{(1-\nu_s)(1-2\nu_s)}{(1+\nu_s)[(1-\nu_s)^2\alpha_n^2-(1-2\nu_s)]}.
\label{1bi(2.9A)}
\end{equation}
Here, $\alpha_n$ are the roots of the characteristic equation (with $J_0$ and $J_1$ being Bessel functions of the first kind)
$$
J_0(x)-\frac{(1-2\nu_s)}{(1-\nu_s)}\frac{J_1(x)}{x}=0.
$$

The inverse relation for Eq.~(\ref{1bi(2.5)}) can be represented as
\begin{equation}
\epsilon(\hat{t})=\int\limits_{0-}^{\hat{t}} \frac{d \hat{F}(\tau)}{d\tau}
\hat{M}(\hat{t}-\tau)\,d\tau,
\label{1bi(2.10)}
\end{equation}
where $\hat{M}(\hat{t})=\mathcal{L}^{-1}\{\bar{M}(s)\}$ with $\bar{M}(s)$ defined by the formula 
\begin{equation}
\bar{M}(s)=\frac{\displaystyle I_0(\sqrt{s})-\frac{2\mu_s}{H_A}\frac{I_1(\sqrt{s})}{\sqrt{s}}
}{\displaystyle s\Bigl(3 I_0(\sqrt{s})-\frac{8\mu_s}{H_A}\frac{I_1(\sqrt{s})}{\sqrt{s}}\Bigr)}.
\label{1bi(2.11)}
\end{equation}

Again, making use of the results by \citet{Armstrong_et_al_1984}, we get
\begin{equation}
\hat{M}(\hat{t})=\frac{1}{2(1+\nu_s)}-\sum_{n=1}^\infty
\frac{B_n}{2(1+\nu_s)}\exp(-\beta_n^2\hat{t}),
\label{1bi(2.12)}
\end{equation}
\begin{equation}
B_n=\frac{4(1-\nu_s^2)(1-2\nu_s)}{9(1-\nu_s)^2\beta_n^2-8(1+\nu_s)(1-2\nu_s)},
\label{1bi(2.12B)}
\end{equation}
where $\beta_n$ are the roots of the corresponding characteristic equation
$$
J_0(x)-\frac{4(1-2\nu_s)}{3(1-\nu_s)}\frac{J_1(x)}{x}=0.
$$

The short-time asymptotic approximation for the kernel $\hat{K}(\hat{t})$ (the same can be done with $\hat{M}(\hat{t})$) can be obtained by evaluating the inverse of $\bar{K}(s)$ as $s\to\infty$. For this purpose, we apply the well known asymptotic formula (see, e.g., \citep{GradshteynRyzhik})
\begin{equation}
I_n(z)=\frac{e^z}{\sqrt{2\pi z}}\Bigl\{1+\frac{(1-4n^2)}{8z}+O(z^{-2})\Bigr\}.
\label{1bi(2.6)}
\end{equation}
Making use of (\ref{1bi(2.6)}), we expand the right-hand sides of (\ref{1bi(2.4)}) and (\ref{1bi(2.11)}) in terms of $1/\sqrt{s}$. As a result, we arrive at the following asymptotic expansions (note a misprint in formula (36b) in \citep{Armstrong_et_al_1984}):
\begin{equation}
\hat{K}(\hat{t})=3-\frac{2(1-2\nu_s)}{\sqrt{\pi}(1-\nu_s)}
\sqrt{\hat{t}}+O(\hat{t}), \quad \hat{t}\to 0,
\label{1bi(2.7)}
\end{equation}
\begin{equation}
\hat{M}(\hat{t})=\frac{1}{3}+\frac{2(1-2\nu_s)}{9\sqrt{\pi}(1-\nu_s)}
\sqrt{\hat{t}}+O(\hat{t}), \quad \hat{t}\to 0.
\label{1bi(2.8)}
\end{equation}

The asymptotic approximations (\ref{1bi(2.7)}) and (\ref{1bi(2.8)}) can be used for evaluating impact unconfined compression test with the impact duration relatively small compared to the so-called gel diffusion time for the biphasic material $t_g=a^2/(H_A k)$.

In the dimensional form, Eqs.~(\ref{1bi(2.5)}) and (\ref{1bi(2.10)}) can be recast as follows:
\begin{equation}
F(t)=\frac{\pi a^2 E_s}{h}\int\limits_{0-}^t \frac{d w(\tau)}{d\tau} K(t-\tau)\,d\tau,
\label{1bi(2.13)}
\end{equation}
\begin{equation}
w(t)=\frac{h}{\pi a^2 E_s}\int\limits_{0-}^t \frac{d F(\tau)}{d\tau} M(t-\tau)\,d\tau.
\label{1bi(2.14)}
\end{equation}
Here we introduced the notation 
\begin{equation}
K(t)=1+\sum_{n=1}^\infty A_n\exp\Bigl(-\frac{t}{\rho_n}\Bigr),
\label{1bi(2.15)}
\end{equation}
\begin{equation}
M(t)=1-\sum_{n=1}^\infty B_n\exp\Bigl(-\frac{t}{\tau_n}\Bigr),
\label{1bi(2.16)}
\end{equation}
\begin{equation}
\rho_n=\frac{a^2}{\alpha_n^2 H_Ak}, \quad
\tau_n=\frac{a^2}{\beta_n^2 H_Ak}.
\label{1bi(2.16RT)}
\end{equation}
Note that the sequences $\rho_1>\rho_2>\ldots $ and $\tau_1>\tau_2>\ldots $, which are defined by formulas (\ref{1bi(2.16RT)}), represent the discrete relaxation and retardation spectra, respectively. Recall that the coefficients $A_n$ and $B_n$ are defined by formulas (\ref{1bi(2.9A)}) and (\ref{1bi(2.12B)}).

Further, in view of (\ref{1bi(2.7)}) and (\ref{1bi(2.8)}), we will have 
\begin{equation}
K_0=K(0)=\frac{3}{2(1+\nu_s)}, 
\label{1bi(2.17K0)}
\end{equation}
\begin{equation}
M_0=M(0)=\frac{2}{3}(1+\nu_s).
\label{1bi(2.17M0)}
\end{equation}
Hence, the following equalities take place:
\begin{equation}
1+\sum_{n=1}^\infty A_n=\frac{3}{2(1+\nu_s)}, 
\label{1bi(2.18An)}
\end{equation}
\begin{equation}
1-\sum_{n=1}^\infty B_n=\frac{2}{3}(1+\nu_s).
\label{1bi(2.18Bn)}
\end{equation}

Now, using formula (\ref{1bi(2.17M0)}), we can rewrite the function (\ref{1bi(2.16)}) as 
\begin{equation}
M(t)=M_0+\sum_{n=1}^\infty B_n\Bigl(1-\exp\Bigl(-\frac{t}{\tau_n}\Bigr)\Bigr).
\label{1bi(2.19)}
\end{equation}

By analogy with the viscoelastic model, the functions $K(t)$ and $M(t)$ will be called the biphasic relaxation and creep functions for unconfined compression.

\section{Cyclic compressive loading}
\label{1biSection3}

Following \citet{Li_et_al_1995}, we assume that a biphasic tissue sample is subjected to a cyclic displacement input
\begin{equation}
w(t)=[w_0(1-\cos\omega t)+w_1]H(t),
\label{1bi(3.1)}
\end{equation}
where $H(t)$ is the Heaviside function, $w_1/h$ is the prestrain resulting from the initial deformation applied to the sample to create the desired preload, $w_0$ is the displacement amplitude, i.e., $w_0/h$ is equal to one-half the peak-to-peak cyclic strain input superimposed on the prestrain, and $\omega=2\pi f$ is the angular frequency with $f$ being the loading frequency. 

Differentiating (\ref{1bi(3.1)}), we get
\begin{equation}
\frac{dw(t)}{dt}=H(t)w_0\omega\sin\omega t+w_1\delta(t).
\label{1bi(3.2)}
\end{equation}
Now, substituting the expression (\ref{1bi(3.2)}) into Eq.~(\ref{1bi(2.13)}), we arrive, after some algebra, at the following resulting stress output:
\begin{eqnarray}
\frac{F(t)}{\pi a^2} & = & \frac{E_s}{h}\Biggl\{ w_1 K(t)
\nonumber\\
{} & {} & {}+w_0\biggl(
1+\sum_{n=1}^\infty \frac{\rho_n^2\omega^2 A_n}{1+\rho_n^2\omega^2}\exp\Bigl(-\frac{t}{\rho_n}\Bigr)
\biggr) \nonumber\\
{} & {} & {}-w_0\bigl[
K_1(\omega)\cos\omega t-K_2(\omega)\sin\omega t
\bigr]\Biggr\}.
\label{1bi(3.3)}
\end{eqnarray}
Here we introduced the notation 
\begin{equation}
K_1(\omega)=1+\sum_{n=1}^\infty \frac{\rho_n^2\omega^2 A_n}{1+\rho_n^2\omega^2},
\label{1bi(3.4K1)}
\end{equation}
\begin{equation}
K_2(\omega)=\sum_{n=1}^\infty \frac{\rho_n\omega A_n}{1+\rho_n^2\omega^2}.
\label{1bi(3.4K2)}
\end{equation}

To assign a physical meaning to the introduced functions $K_1(\omega)$ and $K_2(\omega)$, let us compare the oscillating part of the input strain, that is $-(w_0/h)\cos\omega t$, with the corresponding oscillating part of the compressive stress, which is equal to 
$-E_s(w_0/h)\bigl[K_1(\omega)\cos\omega t-K_2(\omega)\sin\omega t\bigr]$.
By analogy with the viscoelastic model, we obtain that $K_1(\omega)$ and $K_2(\omega)$ represent, respectively, the apparent relative storage and loss moduli. Correspondingly, the apparent loss angle, $\delta(\omega)$, can be introduced by the formula
\begin{equation}
\tan\delta(\omega)=\frac{K_2(\omega)}{K_1(\omega)}. 
\label{1bi(3.5)}
\end{equation}
The apparent loss angle $\delta(\omega)$ describes the phase difference between the displacement input and force output. 

In the case of load-controlled compression, following \citet{Suh_et_al_1995}, we will assume that the tissue sample is subjected to a cyclic compressive loading 
\begin{equation}
F(t)=[F_0(1-\cos\omega t)+F_1]H(t),
\label{1bi(3.6)}
\end{equation}
where $F_0$ is the force amplitude, and $F_1$ is the initial preload. 

After substitution of the expression (\ref{1bi(3.6)}) into Eq.~(\ref{1bi(2.14)}), we finally obtain the following resulting strain output:
\begin{eqnarray}
\frac{w(t)}{h} & = & \frac{1}{\pi a^2 E_s}\Biggl\{ F_1 M(t)
\nonumber\\
{} & {} & {}+F_0\biggl(
1-\sum_{n=1}^\infty \frac{\tau_n^2\omega^2 B_n}{1+\tau_n^2\omega^2}\exp\Bigl(-\frac{t}{\tau_n}\Bigr)
\biggr) \nonumber\\
{} & {} & {}-F_0\bigl[
M_1(\omega)\cos\omega t+M_2(\omega)\sin\omega t
\bigr]\Biggr\}.
\label{1bi(3.7)}
\end{eqnarray}
Here we introduced the notation 
\begin{equation}
M_1(\omega)=1-\sum_{n=1}^\infty \frac{\tau_n^2\omega^2 B_n}{1+\tau_n^2\omega^2},
\label{1bi(3.5M1)}
\end{equation}
\begin{equation}
M_2(\omega)=\sum_{n=1}^\infty \frac{\tau_n\omega B_n}{1+\tau_n^2\omega^2}.
\label{1bi(3.5M2)}
\end{equation}
Note that $M_1(\omega)$ and $M_2(\omega)$ have a physical meaning of the apparent relative storage and loss compliances, respectively. 

\section{Displacement-controlled unconfined compression test}
\label{1biSection4}

Consider an unconfined compression test with the upper plate displacement specified according to the equation 
\begin{equation}
w(t)=w_0\sin\omega t, \quad t\in(0,\pi/\omega).
\label{1bi(4.1)}
\end{equation}
The maximum displacement, $w_0$, will be achieved at the time moment $t_m=\pi/(2\omega)$. The moment of time $t=\tilde{t}_M^\prime$, when the contact force $F(t)$ vanishes, determines the duration of the contact. The contact force itself can be evaluated according to Eqs.~(\ref{1bi(2.13)}) and (\ref{1bi(4.1)}) as follows:
\begin{equation}
F(t)=\frac{\pi a^2 E_s}{h}w_0\omega \int\limits_0^t \cos\omega\tau K(t-\tau)\,d\tau.
\label{1bi(4.2)}
\end{equation}

According to Eq.~(\ref{1bi(4.2)}), the contact forced at the moment of maximum displacement is given by
\begin{equation}
F(t_m)=\frac{\pi a^2 E_s}{h}w_0\tilde{K}_1(\omega),
\label{1bi(4.3)}
\end{equation}
where we introduced the notation
\begin{equation}
\tilde{K}_1(\omega)=\omega \int\limits_0^{\pi/(2\omega)} K(\tau)\sin\omega \tau\,d\tau.
\label{1bi(4.4)}
\end{equation}
By analogy with the viscoelastic case, the quantity $\tilde{K}_1(\omega)$ will be called the incomplete apparent storage modulus. 

Substituting the expression (\ref{1bi(2.15)}) into the right-hand side of Eq.~(\ref{1bi(4.4)}), we obtain
\begin{equation}
\tilde{K}_1(\omega)=1+ \sum_{n=1}^\infty \frac{\omega\rho_n A_n}{\omega^2\rho_n^2+1}\biggl(
\omega\rho_n-\exp\Bigl(-\frac{\pi}{2\omega\rho_n}\Bigr)\biggr).
\label{1bi(4.5)}
\end{equation}

Now, taking into consideration Eqs.~(\ref{1bi(3.4K1)}) and (\ref{1bi(4.5)}), we may conclude that the difference between the apparent storage modulus $K_1(\omega)$ and the incomplete apparent storage modulus $\tilde{K}_1(\omega)$ is relatively small at low frequencies. To be more precise, the difference 
$K_1(\omega)-\tilde{K}_1(\omega)$ is positive and of order $O\bigl(\omega\rho_1\exp(-\pi/(2\omega\rho_1))\bigr)$ as $\omega\to 0$, where $\rho_1$ is the maximum relaxation time. 

In the high frequency limit, the upper limit of the integral (\ref{1bi(4.4)}) tends to zero as $\omega$ increases. Thus, the  behavior of $\tilde{K}_1(\omega)$ as $\omega\to+\infty$ will depend on the behavior of $K(t)$ as $t\to 0$. According to (\ref{1bi(2.7)}), as $\omega\to\infty$, we will have
\begin{equation}
\tilde{K}_1(\omega)=K_0-\frac{s_{1/2}(1-2\nu_s)}{\sqrt{\pi}(1-\nu_s^2)}\sqrt{\frac{H_A k}{a^2}}
\frac{1}{\sqrt{\omega}}+O(\omega^{-1}),
\label{1bi(4.6)}
\end{equation}
where
$
s_{1/2}=\int_0^{\pi/2} \sqrt{x}\sin x\,dx.
$

On the other hand, due to (\ref{1bi(3.4K1)}), (\ref{1bi(2.17K0)}), and (\ref{1bi(2.18An)}), the following limit relation holds true:
$\lim K_1(\omega)=K_0$ as $\omega\to\infty$.
Thus, in view of (\ref{1bi(4.6)}), we conclude that 
$\tilde{K}_1(\omega)\simeq K_1(\omega)$ for $\omega\to\infty$ as well as 
$\tilde{K}_1(\omega)\simeq K_1(\omega)$ for $\omega\to 0$.
In other words, the incomplete apparent storage modulus $\tilde{K}_1(\omega)$ obeys both asymptotic behaviors of the apparent storage modulus $K_1(\omega)$.

\section{Force-controlled unconfined compression test}
\label{1biSection5}

Consider now an unconfined compression test with the external force specified according to the equation 
\begin{equation}
F(t)=F_0\sin\omega t, \quad t\in(0,\pi/\omega).
\label{1bi(5.1)}
\end{equation}
The maximum contact force, $F_0$, will be achieved at the time moment $t_M=\pi/(2\omega)$. The moment of time $t_M^\prime=\pi/\omega$, when the contact force $F(t)$ vanishes, determines the duration of the compression test. According to Eqs.~(\ref{1bi(2.14)}) and (\ref{1bi(5.1)}), the upper plate displacement can be evaluated as follows:
\begin{equation}
w(t)=\frac{h}{\pi a^2 E_s}F_0\omega \int\limits_0^t \cos\omega\tau M(t-\tau)\,d\tau.
\label{1bi(5.2)}
\end{equation}

Due to Eq.~(\ref{1bi(5.2)}), the displacement at the moment of maximum contact force is given by
\begin{equation}
w(t_M)=\frac{h}{\pi a^2 E_s}F_0 \tilde{M}_1(\omega),
\label{1bi(5.3)}
\end{equation}
where we introduced the notation
\begin{equation}
\tilde{M}_1(\omega)=\omega \int\limits_0^{\pi/(2\omega)} M(\tau)\sin\omega \tau\,d\tau.
\label{1bi(5.4)}
\end{equation}
By analogy with the viscoelastic case, the quantity $\tilde{M}_1(\omega)$ will be called the incomplete apparent storage compliance. 

Substituting the expression (\ref{1bi(2.16)}) into the right-hand side of Eq.~(\ref{1bi(5.4)}), we obtain
\begin{equation}
\tilde{M}_1(\omega)=M_0+ \sum_{n=1}^\infty \frac{B_n}{\omega^2\tau_n^2+1}\biggl(
1+\omega\tau_n\exp\Bigl(-\frac{\pi}{2\omega\tau_n}\Bigr)\biggr).
\label{1bi(5.5)}
\end{equation}

In the same way as it was done with $\tilde{K}_1(\omega)$, it can be shown that the incomplete apparent storage compliance $\tilde{M}_1(\omega)$ obeys both asymptotic behaviors of the apparent storage compliance $M_1(\omega)$, that is 
$\tilde{M}_1(\omega)\simeq M_1(\omega)$ for $\omega\to 0$ as well as 
$\tilde{M}_1(\omega)\simeq M_1(\omega)$ for $\omega\to\infty$.


\section{Acknowledgment}

The financial support from the European Union Seventh Framework Programme under contract number PIIF-GA-2009-253055 is gratefully acknowledged. 



\begin{thebibliography}{}

%
\bibitem[Appleyard et al.(2001)]{Appleyard_et_al_2001}
Appleyard, R.C., Swain, M.V., Khanna, S., Murrell, G.A.C.:
The accuracy and reliability of a novel handheld dynamic indentation probe for analysing articular cartilage.
Phys. Med. Biol. {\bf 46}, 541--550 (2001)
%
\bibitem[Argatov(2012)]{Argatov2012}
Argatov, I.: 
Sinusoidally-driven flat-ended indentation of time-dependent materials: Asymptotic models for low and high rate loading. 
Mech. Mater. {\bf 48}, 56--70 (2012)

%
\bibitem[Armstrong et al.(1984)]{Armstrong_et_al_1984}
Armstrong, C.G., Lai, W.M., Mow, V.C.:
An analysis of the unconfined compression of articular cartilage. 
J. Biomech. Eng. {\bf 106}, 165--173 (1984)
%
\bibitem[Butcher and Segalman(1984)]{ButcherSegalman2000}
Butcher, E.A., Segalman, D.J.:
Characterizing damping and restitution in compliant impacts via modified K-V and higher-order linear viscoelastic models.
J. Appl. Mech. {\bf 67}, 831--834 (2000)

\bibitem[Gradshteyn and Ryzhik(1980)]{GradshteynRyzhik}
Gradshteyn, I.S., Ryzhik, I.M.:
Table of Integrals, Series, and Products. Academic, New York (1980)

%
\bibitem[Li et al.(1995)]{Li_et_al_1995}
Li, S., Patwardhan, A.G., Amirouche, F.M.L., Havey, R., Meade, K.P.:
Limitations of the standard linear solid model of intervertebral discs subject to prolonged loading and low-frequency vibration in axial compression.
J. Biomech. {\bf 28}, 779--790 (1995)

%
\bibitem[Mow et al.(1980)]{Mow_et_al_1980}
Mow, V.C., Kuei, S.C., Lai, W.M., Armstrong, C.G.:
Biphasic creep and stress relaxation of articular cartilage in compression: Theory and experiments.
J. Biomech. Eng. {\bf 102}, 73--84 (1980)
%
\bibitem[Stolz et al.(2007)]{Stolz_et_al_2007}
Stolz, M., Aebi, U.,  Stoffler, D.:
Developing scanning probe-based nanodevices --- stepping out of the laboratory into the clinic.
Nanomed.: Nanotechnol. Biol. Med. {\bf 3}, 53--62 (2007)
%
\bibitem[Suh et al.(1995)]{Suh_et_al_1995}
Suh, J.-K., Li, Z., Woo, S.L.-Y.:
Dynamic behavior of a biphasic cartilage model under cyclic compressive loading.
J. Biomech. {\bf 28}, 357--364 (1995)







\end{thebibliography}
\end{document}